\documentclass[11pt]{article}
\usepackage{amsmath}
\usepackage{amsfonts}
\usepackage{mathrsfs}
\usepackage{amsfonts}
\usepackage{amssymb,amsfonts,amsmath, psfrag,eepic,colordvi,epsfig}
\topskip 
\numberwithin{equation}{section}
\newtheorem{theorem}{Theorem}[section]

\newtheorem{conjecture}[theorem]{Conjecture}

\newtheorem{lemma}[theorem]{Lemma}

\begin{document}
\parskip 7pt

\pagenumbering{arabic}
\def\sof{\hfill\rule{2mm}{2mm}}
\def\ls{\leq}
\def\gs{\geq}
\def\SS{\mathcal S}
\def\qq{{\bold q}}
\def\MM{\mathcal M}
\def\TT{\mathcal T}
\def\EE{\mathcal E}
\def\lsp{\mbox{lsp}}
\def\rsp{\mbox{rsp}}
\def\pf{\noindent {\it Proof.} }
\def\mp{\mbox{pyramid}}
\def\mb{\mbox{block}}
\def\mc{\mbox{cross}}
\def\qed{\hfill \rule{4pt}{7pt}}
\def\pf{\noindent {\it Proof.} }
\textheight=22cm

{\Large
\begin{center}
Proofs of two conjectures on
 congruences of overcubic partition triples
\end{center}
}

\begin{center}

  Jiayu Chen$^{1}$, Jing Jin$^{2}$,
  and Olivia X.M. Yao$^{3}$

$^{1,3}$School of Mathematical Sciences, \\
  Suzhou University of Science and
Technology, \\
 Suzhou,  215009, Jiangsu,
 P. R. China

  $^{2}$College of Agricultural
   Information,\\
   Jiangsu Agri-animal
  Husbandry Vocational College,
  \\
  Taizhou, 225300,  Jiangsu,  P. R. China

email:
  $^{1}$ataraxijychen@163.com,
   $^{2}$jinjing19841@126.com,
   $^{3}$yaoxiangmei@163.com

 \end{center}

\noindent {\bf Abstract.}
  Let $\overline{bt}(n)$
   denote
    the number of overcubic partition
triples of $n$. Nayaka, Dharmendra
 and Kumar proved some congruences
  modulo 8, 16 and 32 for $\overline{bt}(n)$.
  Recently, Saikia and Sarma
   established  some congruences
    modulo 64
     for $\overline{bt}(n)$ by using both
      elementary techniques and
      the theory of modular forms.
    In their paper, they also  posed two
 conjectures on infinite families
  of  congruences modulo 64
  and 128 for $\overline{bt}(n)$.
   In this paper, we confirm the
   two conjectures.



\section{Introduction}

\allowdisplaybreaks

Given a positive integer
 $n$,
  a partition
     of  $n$  is a finite
      weakly decreasing sequence
      of positive integers
$\pi=(\pi_{1},\pi_{2},\ldots,\pi_{k})$ such that
$\pi_{1}+\pi_{2}+\cdots+\pi_{k}=n$.
 As usual, let $p(n)$  denote
  the number of partitions of
    $n$ and set $p(0)=1$. Euler found that the generating
       function for $p(n)$
        is
\begin{align*}
\sum_{n=0}^\infty p(n) q^n=\frac{1}{f_1}.
\end{align*}
Here and  throughout the paper, we use the
 following notation:
\begin{align*}
f_k:=\prod_{n=1}^\infty (1-q^{nk}).
\end{align*}

 Recall that
   the partitions in which even parts
    come in   two colors blue (denoted by $b$)
    and red  (denoted by $r$)
   are known as cubic partitions.
    For instance, the nine cubic
     partitions of 4 are:
\[
4_b, \quad 4_r,\quad 3+1, \quad
 2_b+2_b, \quad 2_b+2_r, \quad
 2_r+2_r, \quad 2_r+1+1, \quad
 2_b+1+1,\quad 1+1+1+1.
\]
Let $a(n)$ denote the number
 of cubic partitions of $n$.
  The generating
 function for $a(n)$ is
\[
  \sum_{n=0}^\infty
     a(n)q^n =\frac{1}{f_1f_2
      }.
\]

In a series of papers, Chan \cite{Chan-1,Chan-2,Chan-3}
 studied
congruence properties for $a(n)$ and
   proved some congruences modulo
     powers of 3 for $a(n)$. Zhao and Zhong
     \cite{Zhao}
       proved some congruences
        on cubic partition pairs.
     In 2010,  Kim
 \cite{Kim-2010} studied the overcubic
       partition function    $\overline{a}
      (n)$  which
  counts all of the overlined version
  of the cubic partitions
counted by $a(n)$, namely,
    the cubic partitions where the
first instance of each part is
 allowed to be overlined. The generating function
  of $\overline{a}(n)$ is
\[
\sum_{n=0}^\infty
 \overline{a}(n)q^n=\frac{f_4}{f_1^2f_2}.
\]
Based on the theory
 of modular forms,
  Kim proved that
\[
\sum_{n=0}^\infty
 \overline{a}(3n+2)q^n=3\frac{f_3^6f_4^3
  }{f_1^8f_2^3}.
\]
 A number of
  congruences for $\overline{a}(n)
 $
  have been proved by Sellers \cite{Sellers}.
 In 2012, Kim \cite{Kim-2012} investigated
  congruence properties of $\overline{b}(n)
   $, which counts  the
number of overcubic partition pairs of $n$.
 Note that the generating function
  of $ \overline{b}(n)
   $ is
\[
\sum_{n=0}^\infty
 \overline{b}(n)q^n=\frac{f_4^2}{f_1^4f_2^2}.
\]

Recently, Nayaka, Dharmendra and Kumar
 \cite{Nayaka} investigated  congruence
  properties  for $\overline{bt}(n)
   $, which counts  the
number of overcubic partition triples
  of $n$. The generating function of
   $\overline{bt}(n)
   $ is    \begin{align}
\sum_{n=0}^\infty
 \overline{bt}(n)q^n=\frac{f_4^
 3}{f_1^6f_2^3}. \label{1-1}
  \end{align}
  They also proved a number of congruences
   modulo 8, 16 and 32 for $\overline{bt}(n)$.
    For example, they proved that for $n\geq 0
     $,
     \begin{align*}
\overline{bt}(8n+5) \equiv &\;
 0\pmod 8, \\
 \overline{bt}(16n+10) \equiv &\;
 0\pmod {16}, \\
 \overline{bt}(8n+7) \equiv &\;
 0\pmod {32}.
     \end{align*}
Very recently, Saikia and Sarma \cite{Saikia}
 proved  many  congruences modulo
  64 and 128 for $\overline{bt}(n)$
  by using both elementary techniques
   and the theory of modular forms.
    For example, they proved that
     for $n\geq 0$,
     \begin{align}
\overline{bt}(8n+7) \equiv &\;
 0\pmod {64}, \label{v-1}\\
 \overline{bt}(16n+14) \equiv &\;
 0\pmod {64}, \label{v-2}\\
 \overline{bt}(32n+28) \equiv &\;
 0\pmod {64}, \label{v-3}\\
 \overline{bt}(72n+21) \equiv &\;
 0\pmod {128}, \label{v-4}\\
 \overline{bt}(72n+69) \equiv &\;
 0\pmod {384}. \label{v-5}
     \end{align}
In their paper, they posed
 the following two conjectures
  on infinite families
   of congruences modulo 64 and 128.

\begin{conjecture}\label{c-1}
 \cite{Saikia} For $n,\alpha\geq 0$,
 \begin{align}
 \overline{bt}(2^{\alpha}(8n+7)) \equiv &\;
 0\pmod {64}.\label{1-2}
 \end{align}
\end{conjecture}

\begin{conjecture}\label{c-2}
  \cite{Saikia} For $n,\alpha\geq 0$,
 \begin{align}
 \overline{bt}(144n+42) \equiv &\;
 0\pmod {384},\label{1-3}\\
  \overline{bt}(2^{\alpha}(72n+21)) \equiv &\;
 0\pmod {128},\label{1-4}\\
   \overline{bt}(2^{\alpha}(72n+69)) \equiv &\;
 0\pmod {128}. \label{1-5}
 \end{align}
\end{conjecture}

The aim of this paper is to confirm
 Conjectures
 \ref{c-1} and \ref{c-2}.

\section{Preliminaries}

To prove Conjectures \ref{c-1}
 and \ref{c-2}, we first prove
 the following lemma.

\begin{lemma}
We have
\begin{align}
  \sum_{n=0}^\infty
 \overline{bt}(4n)q^n&\; \equiv
 32q\frac{f_2^{47}}{
 f_4f_8^2} \cdot
 \left(\frac{1}{f_1^{4}}\right)^{12}
 \cdot \frac{1}{f_1^2}+56q\frac{
 f_2^{61}f_8^2}{ f_4^{15}} \cdot
 \left(\frac{1}{f_1^{4}}\right)^{13}
 \cdot \frac{1}{f_1^2}
 \nonumber\\[6pt]
&\;\quad   +\frac{f_2^{71}}{ f_4^{17}
   f_8^2} \cdot
 \left(\frac{1}{f_1^{4}}\right)^{14}
 \cdot \frac{1}{f_1^2} \pmod {128},
  \label{2-1} \\[6pt]
   \sum_{n=0}^\infty
 \overline{bt}(8n)q^n&\; \equiv
 32q\frac{f_2^{169}f_8^2
  }{
 f_4^{51}} \cdot
 \left(\frac{1}{f_1^{4}}\right)^{31}
 \cdot \frac{1}{f_1^2}+16q\frac{
 f_2^{155} }{f_4^{37}f_8^2} \cdot
 \left(\frac{1}{f_1^{4}}\right)^{30}
 \cdot \frac{1}{f_1^2}
  \nonumber\\[6pt]
&\;\quad+\frac{f_2^{179}}{f_4^{53}
   f_8^2} \cdot
 \left(\frac{1}{f_1^{4}}\right)^{32}
 \cdot \frac{1}{f_1^2} \pmod {128},
  \label{2-2} \\[6pt]
  \sum_{n=0}^\infty
   \overline{bt}(16n)q^n&\; \equiv
 96q\frac{f_2^{385}f_8^2
  }{
 f_4^{123}}\cdot
 \left(\frac{1}{f_1^{4}}\right)^{67}
 \cdot \frac{1}{f_1^2}
  +\frac{f_2^{395}}{ f_4^{125}
   f_8^2} \cdot
 \left(\frac{1}{f_1^{4}}\right)^{68}
 \cdot \frac{1}{f_1^2}\pmod {128},
  \label{2-3}\\[6pt]
  \sum_{n=0}^\infty
   \overline{bt}(32n)q^n&\; \equiv
 96q\frac{f_2^{803}
  }{
 f_4^{253}f_8^2} \cdot
 \left(\frac{1}{f_1^{4}}\right)^{138}
 \cdot \frac{1}{f_1^2}
  +96q\frac{
  f_2^{817}f_8^2}{
     f_4^{267}
    }\cdot
 \left(\frac{1}{f_1^{4}}\right)^{139}
 \cdot \frac{1}{f_1^2}
 \nonumber\\[6pt]
&\;\quad+\frac{f_2^{827}}{ f_4^{269}
   f_8^2}  \cdot
 \left(\frac{1}{f_1^{4}}\right)^{140}
 \cdot \frac{1}{f_1^2} \pmod {128}.
 \label{2-4}
\end{align}
\end{lemma}

\noindent{\it Proof.} We can rewrite
 \eqref{1-1} as
 \begin{align}\label{2-5}
 \sum_{n=0}^\infty
 \overline{bt}(n)q^n=\frac{f_4^3}{f_2^3}
 \cdot
 \frac{1}{f_1^4}\cdot \frac{1}{f_1^2}.
 \end{align}
 It follows from \cite[Entry 25, pp. 40]{Berndt-1} that
\begin{align}
 \varphi(q)+
  \varphi(-q)&=2\varphi(q^4)
 \label{3-10}
,\\
 \varphi(q)-\varphi(-q)&=4q\psi(q^8).
 \label{3-11} \\
 \varphi(q)^2-\varphi(-q)^2  &= 8q\psi(q^4)^2,
 \label{4-5}\\
 \varphi(q)^2+\varphi(-q)^2
  &= 2\varphi(q^2)^2
 ,\label{a-1}
\end{align}
where
\begin{align}
\varphi(q):&=\sum_{n=-\infty}^\infty q^{n^2}
 =\frac{f_2^5}{f_1^2f_4^2}, \label{3-3}
 \\
  \psi(q):&=\sum_{n=0}^\infty q^{n(n+1)/2}
 =\frac{f_2^2}{f_1}. \label{3-4}
\end{align}
In view of \eqref{3-10}--\eqref{a-1},
 \begin{align}
 \frac{1}{f_1^2}=&\;
\frac{f_8^5}{f_2^5f_{16}^2} +2q\frac{f_4^2f_{16}^2}{f_2^5f_8}
\label{2-6}
\end{align}
 and
 \begin{align} \frac{1}{f_1^4}=&\;
\frac{f_4^{14}}{f_2^{14}f_8^4}
+4q\frac{f_4^2f_8^4}{f_2^{10}}.
\label{2-7}
\end{align}
Substituting \eqref{2-6}
 and \eqref{2-7} into \eqref{2-5} and  extracting those
terms in which the power of  $q
 $ is congruent to 0 modulo 2, then
 replacing $q^{2}$ by $q$, we arrive at
 \begin{align}
 \sum_{n=0}^\infty
 \overline{bt}(2n)q^n&\;
 =8q f_2^7
  f_4^3 f_8^2 \cdot
   \left(
    \frac{1}{f_1^4}\right)^4
    \cdot \frac{1}{f_1^2} +\frac{
  f_2^{17} f_4}{f_8^2}\left(
    \frac{1}{f_1^4}\right)^5
    \cdot \frac{1}{f_1^2}. \label{2-8}
  \end{align}
Substituting  \eqref{2-6}
 and \eqref{2-7} into \eqref{2-8}
  and exacting
 those terms in which the power
  of $q$ is congruent to 0
  modulo 2,
 then  replacing
  $q^2
 $ by $q$,
 we arrive at \eqref{2-1}.
  Substituting \eqref{2-6} and
  \eqref{2-7} into \eqref{2-1}
 and picking out
 those terms in which the power
  of $q$ is congruent to 0
  modulo 2,
 then  replacing
  $q^2
 $ by $q$,
  we arrive at \eqref{2-2}.
 Substituting \eqref{2-6} and \eqref{2-7}
   into \eqref{2-2}
 and exacting
 those terms in which the power
  of $q$ is congruent to 0
  modulo 2,
 then  replacing
  $q^2
 $ by $q$,
  we arrive at \eqref{2-3}.
Substituting \eqref{2-6} and \eqref{2-7}
  into \eqref{2-3}
 and picking out
 those terms in which the power
  of $q$ is congruent to 0
  modulo 2,
 then  replacing
  $q^2
 $ by $q$,
 we arrive at \eqref{2-4}.
  This completes the proof of this lemma.
  \qed

 \section{Proof of Conjecture 1.1}

 Using the binomial theorem,
 one can prove that if
  $p$ is  a prime, then for all positive integers
 $m$ and $k$,
\begin{align}
 f_m^{p^k}\equiv f_{mp}^{p^{k-1}}
  \pmod {p^k}.\label{2-9}
\end{align}
In view of  \eqref{2-2}, \eqref{2-3}
  and \eqref{2-9},
\begin{align}
\sum_{n=0}^\infty (
 \overline{bt}(16n)
  -\overline{bt}(8n))q^n\equiv &\;
  -16q\frac{f_2^{15}}{f_1^2f_4}+\frac{
  f_2^{11}f_4^3}{f_1^{18}f_8^2}-\frac{
  f_4^{11}}{f_1^2f_2^{13}f_8^2}
  \nonumber\\[6pt]
  \equiv&\; \frac{f_2^{3}f_4^3
   }{f_1^2 f_8^2
    }\left(-16q \frac{
     f_8^4}{f_4^2}
    +\frac{f_2^8}{f_1^{16}}
     -\frac{f_4^{8}}{f_2^{16}}\right)
      \pmod {64}. \label{3-1}
\end{align}
It is easy to check that
\begin{align}\label{3-2}
 \sum_{m,n=1}^\infty (-1)^{m+n}
 q^{m^2+n^2}=&\;\sum_{m,n=1,\atop
 m>n}^\infty
  (-1)^{m+n}q^{m^2+n^2}
  +\sum_{m,n=1,\atop
 n>m}^\infty
  (-1)^{m+n}q^{m^2+n^2}+\sum_{n=1}^\infty
   q^{2n^2}
   \nonumber\\[6pt]
   =&\;2\sum_{m,n=1,\atop
 m>n}^\infty
  (-1)^{m+n}q^{m^2+n^2}
  +\sum_{n=1}^\infty
   q^{2n^2}.
\end{align}

Replacing $q$ by $-q$ in \eqref{3-3} yields
\begin{align}
\varphi(-q)=\frac{f_1^2}{f_2}=1+2\sum_{n=1}^\infty
 (-1)^n q^{n^2}.\label{3-5}
\end{align}
Combining \eqref{3-2}
 and \eqref{3-5} yields
\begin{align}
 \frac{f_2}{f_1^2}
 =&\;\frac{1}{1+2\sum_{n=1}^\infty
  (-1)^n q^{n^2}}
 \nonumber\\[6pt]
 =&\; 1+\sum_{j=1}^\infty
  (-2)^j \left(\sum_{t=1}^\infty
  (-1)^t q^{t^2}\right)^j
 \nonumber\\[6pt]
 \equiv&\; 1-2\sum_{n=1}^\infty
  (-1)^n q^{n^2}+4\sum_{m,n=1}^\infty
   (-1)^{m+n}q^{m^2+n^2}
 \nonumber\\[6pt]
 \equiv&\; 1-2\sum_{n=1}^\infty
  (-1)^n q^{n^2}+4\sum_{ n=1}^\infty
    q^{2n^2}  \pmod 8. \label{3-6}
\end{align}
In view of \eqref{3-5}
 and  \eqref{3-6},
\begin{align}
\frac{f_2^8}{f_1^{16}}\equiv &\;
 1+48\sum_{n=1}^\infty
 (-1)^n q^{n^2}+32\sum_{n=1}^\infty
  q^{2n^2}+48\left(\sum_{n=1}^\infty
 (-1)^n q^{n^2}\right)^2+32\left(\sum_{n=1}^\infty
 (-1)^n q^{n^2}\right)^4
 \nonumber \\[6pt]
 \equiv &\;
 1+48\sum_{n=1}^\infty
 (-1)^n q^{n^2}+32\sum_{n=1}^\infty
    q^{2n^2}+48\left(\sum_{n=1}^\infty
 (-1)^n q^{n^2}\right)^2
 +32 \sum_{n=1}^\infty
  q^{4n^2}
  \nonumber\\[6pt]
 \equiv &\;
 1+24(\varphi(-q)-1) +16(\varphi( q^2) -1
  ) +12 (\varphi(-q)-1)^2
  +16(\varphi( q^4)-1) \nonumber\\
  \equiv &\; 21 +12\varphi(-q)^2
 +16\varphi( q^2)+16\varphi( q^4)
 \pmod {64},\label{3-7}
\end{align}
where $\varphi(q)$ is defined by \eqref{3-3}.
  Replacing $q$ by $q^2$ in
\eqref{3-7} yields
\begin{align}
\frac{f_4^8}{f_2^{16}}
 \equiv &\;
 21 +12\varphi(-q^2)^2
 +16\varphi( q^4)+16\varphi( q^8)
  \pmod {64}.\label{3-8}
\end{align}
Thanks to \eqref{3-4}, \eqref{3-7}
 and \eqref{3-8},
\begin{align}
-16q \frac{
     f_8^4}{f_4^2}
    +\frac{f_2^8}{f_1^{16}}
     -\frac{f_4^{8}}{f_2^{16}}&\;\equiv -16q\psi(q^4)^2
     +16\varphi( q^2)-16\varphi( q^8)+12\varphi(-q )^2
     -12\varphi(-q^2)^2
     \pmod {64}.\label{3-9}
\end{align}
  In view of \eqref{3-10}
  and \eqref{3-11},
\begin{align}
 & -16q\psi(q^4)^2
       +12\varphi(-q )^2
     -12\varphi(-q^2)^2 +16(\varphi( q^2)-\varphi( q^8))
      \nonumber\\
  =& -16q\psi(q^4)^2
      +12\varphi(-q)
    (\varphi(-q) -\varphi(q))
    +16(\varphi(q^2)- \varphi(q^8))
     \nonumber\\[6pt]
      =& -16q\psi(q^4)^2
      -48q\varphi(-q)\psi(q^8)
     +32q^2\psi(q^{16})
     \nonumber\\[6pt]
    =& -16q\psi(q^4)^2 -48
        q(\varphi(q^4)-2q\psi(q^8)
      )\psi(q^8)
     +32q^2\psi(q^{16})
     \nonumber\\[6pt]
     \equiv& 0 \pmod {64}. \label{3-12}
\end{align}
Here we have used
 the following two results:
 \[
\psi(q^4)^2=  \varphi(q^4)\psi(q^8), \qquad
 \psi(q^8)^2\equiv \psi(q^{16}) \pmod 2.
 \]
In light of \eqref{3-1}, \eqref{3-9} and \eqref{3-12},
\begin{align}
\overline{bt}(16n)
  \equiv \overline{bt}(8n) \pmod {64}.
  \label{3-13}
\end{align}
By \eqref{3-13} and mathematical induction,
 we  deduce that
  for $n,\alpha\geq 0$,
\begin{align}
\overline{bt}(2^{\alpha+3}n)
  \equiv \overline{bt}(8n) \pmod {64}.
  \label{3-14}
\end{align}
Substituting \eqref{2-6} and \eqref{2-7}
   into \eqref{2-2}
 and exacting
 those terms in which the power
  of $q$ is congruent to 1
  modulo 2,
 then dividing by $q$ and  replacing
  $q^2
 $ by $q$,
  we get
  \begin{align}\label{3-15}
\sum_{n=0}^\infty \overline{bt}(16n+8)q^n
 \equiv 48 \frac{f_2^{383}}{
   f_4^{117}f_8^2}\cdot
  \left(\frac{1}{f_1^4
   }\right)^{67}\cdot \frac{1}{f_1^2}
  +2\frac{f_2^{397}f_8^2}{
   f_4^{131}}\cdot
  \left(\frac{1}{f_1^4
   }\right)^{68}\cdot \frac{1}{f_1^2}
   \pmod {128}.
\end{align}
Here we have used    \eqref{2-9}. If we
 substitute   \eqref{2-6} and
\eqref{2-7}
   into \eqref{3-15}
 and pick out
 those terms in which the power
  of $q$ is congruent to 1
  modulo 2,
 then divide by $q$ and
 replace
  $q^2
 $ by $q$,
  we arrive at
\begin{align*}
\sum_{n=0}^\infty \overline{bt}(32n+24)q^n
 &\; \equiv 36 \frac{f_2^{823}f_8^2}{
  f_1^{560}f_4^{271} }
  +32\frac{f_2^{809} }{
  f_1^{556}f_4^{257}f_8^2} \nonumber\\[6pt]
  &\;\equiv  4 \frac{f_4^2}{f_2}
  \cdot \frac{ f_8^2}{
  f_4  }
      \pmod {64}, \qquad ({\rm by}\ \eqref{2-9})
\end{align*}
which yields
\begin{align}\label{3-16}
\overline{bt}(8(8n+7))\equiv 0 \pmod {64}.
\end{align}
Replacing $n$ by $8n+7$
 in \eqref{3-14} and using
  \eqref{3-16} yields
\begin{align}\label{3-17}
\overline{bt}(2^{\alpha+3}(8n+7))\equiv 0 \pmod {64}.
\end{align}
Congruence \eqref{1-2}
  follows from  \eqref{v-1}--\eqref{v-3}
    and \eqref{3-17}.
 This completes
  the proof of  Conjecture \ref{c-1}. \qed

\section{Proof of Conjecture 1.2}

\begin{lemma}\label{L-4-1}
 For $n,\alpha\geq 0$,
 \begin{align}\label{4-1}
\overline{bt}(2^{\alpha+4}n)&\; \equiv \overline{bt}(16n) \pmod
{128}.
 \end{align}

\end{lemma}

\noindent{\it Proof.} In light of
  \eqref{2-3},
  \eqref{2-4} and \eqref{2-9},
\begin{align}\label{4-2}
\sum_{n=0}^\infty (
 \overline{bt}(32n)
  -\overline{bt}(16n))q^n\equiv &\;
 \frac{f_2^{43}
   }{f_1^{18}f_4^{13} f_8^2}  \left(96q \frac{
     f_8^4}{f_4^2}
   +\frac{f_2^{16}}{f_1^{32}}
     -\frac{f_4^{16}}{f_2^{32}}\right)
      \pmod {128}.
\end{align}
By \eqref{3-7},
\begin{align}\label{4-3}
\frac{f_2^{16}}{f_1^{32}}
 \equiv 57+120\varphi(-q)^2+32\varphi(q^2)
 +32\varphi(q^4)+16\varphi(-q)^4 \pmod
 {128}.
\end{align}
Replacing $q$ by $q^2$ in \eqref{4-3} yields
\begin{align}\label{4-4}
\frac{f_4^{16}}{f_2^{32}}
 \equiv 57+120\varphi(-q^2)^2+32\varphi(q^4)
 +32\varphi(q^8)+16\varphi(-q^2)^4
 \pmod {128}.
\end{align}
It is easy to check that
 \begin{align}
\psi(q)^2&=\varphi(q)\psi(q^2).\label{4-6}
 \end{align}
   In light of \eqref{3-10}--\eqref{4-5} and
  \eqref{4-3}--\eqref{4-6},
\begin{align}
&96q \frac{
     f_8^4}{f_4^2}
    +\frac{f_2^{16}}{f_1^{32}}
     -\frac{f_4^{16}}{f_2^{32}}
     \nonumber\\
      =&\;96q\psi(q^4)^2
    +120(\varphi(-q)^2-\varphi(-q^2)^2
     )+32(\varphi(q^2)-\varphi(q^8))+16
     (\varphi(-q)^4-\varphi(-q^2)^4
     )
     \nonumber\\
     =&\;96q\psi(q^4)^2
    +120\varphi(-q)(\varphi(-q)
     -\varphi(q)
     )+32(\varphi(q^2)-\varphi(q^8))+16
     \varphi(-q)^2
      (\varphi(-q)^2-\varphi(q)^2
     )
     \nonumber\\
      =&\;96q\psi(q^4)^2
    -480q \varphi(-q)\psi(q^8)
     +64q^2\psi(q^{16})-128
     \varphi(-q)^2 \psi(q^4)^2
     \nonumber\\
       =&\;96q\psi(q^4)^2
    -480q (\varphi(q^4)-2q\psi(q^8))\psi(q^8)
     +64q^2\psi(q^{16})-128
     \varphi(-q)^2 \psi(q^4)^2
     \nonumber\\
       =&\;96q\psi(q^4)^2
    -480q  \varphi(q^4)\psi(q^8)+960q^2
     \psi(q^8)^2
     +64q^2\psi(q^{16})-128
     \varphi(-q)^2 \psi(q^4)^2
     \nonumber\\
     \equiv&\; 0\pmod {128}. \label{4-7}
\end{align}
Here we have used
 the following congruence:
 \begin{align*}
  \qquad \psi(q^8)^2\equiv
 \psi(q^{16}) \pmod 2.
 \qquad ({\rm by} \eqref{2-9})
 \end{align*}
In view of \eqref{4-2} and \eqref{4-7},
\begin{align}\label{4-8}
\overline{bt}(16n)
  \equiv \overline{bt}(32n) \pmod {128}.
\end{align}
By \eqref{4-8} and  mathematical induction,
 we arrive at \eqref{4-1}. This completes
  the proof of Lemma \ref{L-4-1}. \qed

Now, we are ready to prove Conjecture
 \ref{c-2}.

Substituting \eqref{2-6} and \eqref{2-7}
   into \eqref{2-8}
 and exacting
 those terms in which the power
  of $q$ is congruent to 1
  modulo 2,
 then dividing by $q$ and  replacing
  $q^2
 $ by $q$,
  we arrive at
\begin{align}\label{4-9}
\sum_{n=0}^\infty \overline{bt}(4n+2)q^n
  \equiv &\;  64q \left(
  \frac{1}{f_1^4}\right)^{12}\cdot
  \frac{1}{f_1^2}\cdot
   \frac{f_2^{49}
  f_8^2}{ f_4^7}+28\left(
   \frac{1}{f_1^4}\right)^{13}
   \cdot
   \frac{1}{f_1^2}\cdot
    \frac{f_2^{59}}{
  f_4^9f_8^2}\nonumber\\
  &+2\left(\frac{1}{f_1^4}\right)^{14}\cdot
  \frac{1}{f_1^2}\cdot \frac{f_2^{73}
   f_8^2}{ f_4^{23}} \pmod {128}.
\end{align}
If we substitute \eqref{2-6} and \eqref{2-7}
   into \eqref{4-9}
 and
extract those terms in which
 the power of $q$ is congruent
 to 0 mod 2, then
 replace $q^2$ by $q$, we have
\begin{align}\label{4-10}
\sum_{n=0}^\infty \overline{bt}(8n+2)q^n
  \equiv&\;  64q\cdot
 \left(\frac{1}{f_1^4}\right)^{31}
 \cdot
  \frac{f_2^{163}
  f_8^2}{ f_4^{47} }+96q \cdot
  \left(\frac{1}{f_1^4}\right)^{30}\cdot
   \frac{f_2^{149}}{
 f_4^{33}f_8^2}\nonumber\\
 &+30\left(
 \frac{1}{f_1^4}\right)^{32}\cdot
  \frac{f_2^{173}
   }{ f_4^{49}f_8^2} \pmod {128}.
\end{align}
Substituting \eqref{2-7}
   into \eqref{4-10}
 and exacting
 those terms in which the power
  of $q$ is congruent to 1
  modulo 2,
 then dividing by $q$ and  replacing
  $q^2
 $ by $q$,
  we arrive at
\begin{align}\label{4-11}
\sum_{n=0}^\infty \overline{bt}(16n+10)q^n
 &\; \equiv  32\frac{f_2^{387}}{f_1^{271}f_4^{122}}
  \nonumber\\[6pt]
  &\;\equiv 32\left(\frac{f_2^2}{f_1}\right)^3
   \cdot \frac{f_4^2}{f_2} \pmod {128}.
   \qquad ({\rm by} \eqref{2-9})
\end{align}
It follows from \cite[Corollary (ii), pp. 49]{Berndt-1} that
\begin{align}\label{4-12}
 \frac{f_2^2}{f_1}=&\;
\frac{f_6f_9^2}{f_3f_{18}}+q\frac{f_{18}^2}{f_9}.
\end{align}
 Substituting \eqref{4-12} into
 \eqref{4-11} and picking
out those terms
 in which the power of $q$ is congruent to 2 mod 3, then dividing
 by $q^2$ and
 replacing $q^3$ by $q$ ,we deduce that
\begin{align}\label{4-13}
\sum_{n=0}^\infty \overline{bt}(48n+42)q^n
 &\; \equiv  32q\frac{
 f_6^5f_{12}^2}{f_3^3}
 +32\frac{f_2^3f_3^6f_{12}^2}{f_1^3
 f_6^4}+96\frac{
  f_6^5}{f_{12}} \cdot
  \frac{f_4}{f_1} \nonumber\\
   &\; \equiv  32q\frac{
 f_6^5f_{12}^2}{f_3^3}
 +32\cdot f_1f_2\cdot f_3^6+96\frac{
  f_6^5}{f_{12}} \cdot
  \frac{f_4}{f_1} \pmod {128}.
   \qquad ({\rm by} \eqref{2-9})
\end{align}
In \cite{Andrews}, Andrews,  Hirschhorn
 and Sellers   proved that
 \begin{align}
\frac{f_4}{f_1} =\frac{f_{12}f_{18}^4}
{f_3^3f_{36}^2}
+q\frac{f_6^2f_9^3f_{36}} {f_3^4f_{18}^2}
  +2q^2\frac{f_6f_{18}f_{36}}{f_3^3}.
  \label{4-15}
 \end{align}
 In addition, Hirschhorn
 and Sellers \cite{Hirschhorn}
  proved that
 \begin{align}
f_1f_2=\frac{f_6f_9^4}{f_3f_{18}^2} -qf_9f_{18}
-2q^2\frac{f_3f_{18}^4}{f_6f_9^2}. \label{4-14}
\end{align}
If we substitute
 \eqref{4-15}
  and \eqref{4-14} into \eqref{4-13}
   and pick out  those terms in which
 the power of $q$ is congruent
 to 0 mod 3, then
 replace $q^3$ by $q$, we obtain
\begin{align*}
\sum_{n=0}^\infty \overline{bt}(144n+42)q^n
 &\; \equiv  96 \frac{
 f_2^5 f_6^4}{f_1^3f_{12}^2}
 +32\frac{f_1^5f_2f_3^4}{f_6^2}
 \nonumber\\[6pt]
 &\;\equiv 0 \pmod {128},
 \qquad ({\rm by} \eqref{2-9})
\end{align*}
which yields
\begin{align}
\overline{bt}(144n+42)\equiv 0 \pmod
 {128}. \label{4-16}
\end{align}
 Thanks to \eqref{1-1} and \eqref{2-9},
\begin{align*}
\sum_{n=0}^\infty
 \overline{bt}(n)q^n
&\;\equiv \frac{f_{12}}{f_3^2f_6} \pmod {3},
\end{align*}
which yields
\begin{align}\label{4-18}
\sum_{n=0}^\infty
 \overline{bt}(3n)q^n
  \equiv&\; \frac{f_4}{f_1^2f_2}\nonumber\\
 \equiv &\;\frac{f_1f_4}{
f_2f_3} \pmod {3}.
 \qquad ({\rm by} \eqref{2-9})
\end{align}
Replacing $q$ by $-q$  in \eqref{4-18} yields
\begin{align}\label{4-19}
\sum_{n=0}^\infty
 \overline{bt}(3n) (-1)^n q^n
  \equiv \; \frac{f_2^2}{f_1}
  \cdot
  \frac{f_3f_{12}}{f_6^3} \pmod {3}.
\end{align}
It follows from \eqref{4-12}
 and \eqref{4-19} that for $n\geq 0$,
\begin{align*}
  \overline{bt}(9n+6)
 &\;\equiv 0 \pmod {3}.
\end{align*}
Congruence \eqref{1-3}
 follows from
  \eqref{4-16} and the above congruence.

Substituting \eqref{2-6} and \eqref{2-7}
   into \eqref{2-1}
 and exacting
 those terms in which the power
  of $q$ is congruent to 1
  modulo 2,
 then dividing by $q$ and  replacing
  $q^2
 $ by $q$,
  we arrive at
\begin{align}\label{4-20}
\sum_{n=0}^\infty \overline{bt}(8n+4)q^n
   \equiv  &\; 32q \left(
  \frac{1}{f_1^4}\right)^{30}
  \cdot
  \frac{1}{f_1^2}\cdot
   \frac{f_2^{157}f_8^2
  }{ f_4^{43}}
  +16 \left(\frac{1}{f_1^4}\right)^{31}
  \cdot
  \frac{1}{f_1^2}\cdot
   \frac{f_2^{167}}{
    f_4^{45}f_8^2}
    \nonumber\\
    &\;
    +2\left(\frac{1}{f_1^4}\right)^{32}
  \cdot
  \frac{1}{f_1^2}\cdot\frac{f_2^{181}
    f_8^2}{ f_4^{59}} \pmod {128}.
\end{align}
If we substitute
 \eqref{2-6}
  and \eqref{2-7} into \eqref{4-20}
   and pick out  those terms in which
 the power of $q$ is congruent
 to 0 mod 2, then
 replace $q^2$ by $q$, we get
\begin{align} \label{4-21}
\sum_{n=0}^\infty \overline{bt}(16n+4)q^n
 &\; \equiv  64q \cdot \left(\frac{1}{f_1^4}
 \right)^{67}
 \frac{f_2^{379}f_8^2}{
  f_4^{119}}+18
  \left(\frac{1}{f_1^4}\right)^{68}
  \cdot
   \frac{f_2^{389}}{
  f_4^{121}f_8^2}   \pmod {128}.
\end{align}
Substituting   \eqref{2-7}
   into \eqref{4-21}
 and exacting
 those terms in which the power
  of $q$ is congruent to 1
  modulo 2,
 then dividing by $q$ and  replacing
  $q^2
 $ by $q$,
  we arrive at
\begin{align}
\sum_{n=0}^\infty \overline{bt}(32n+20)q^n
 &\; \equiv 96
 \frac{f_2^{819}}{f_1^{559}
 f_4^{266}}  \pmod {128}. \label{4-22}
\end{align}
Combining \eqref{2-9},
 \eqref{4-11} and \eqref{4-22} yields
 \begin{align}
\overline{bt}(32n+20)\equiv- \overline{bt}(16n+10)\pmod
{128}.\label{4-23}
\end{align}
Replacing $n$ by $9n+2$ in \eqref{4-23}
 and using \eqref{4-16} yields
\begin{align}
\overline{bt}(288n+84)\equiv  0\pmod {128}.
\label{4-24}
\end{align}

Substituting \eqref{2-6} and \eqref{2-7}
   into \eqref{3-15}
 and exacting
 those terms in which the power
  of $q$ is congruent to 0
  modulo 2,
 then   replacing
  $q^2
 $ by $q$,
  we deduce that
\begin{align}
\sum_{n=0}^\infty \overline{bt}(32n+8)q^n
 &\; \equiv  64q\left(
  \frac{1}{f_1^4}\right)^{138}
  \cdot \frac{f_2^{797}}{
 f_4^{249}f_8^2
  }+64q\left(
  \frac{1}{f_1^4}\right)^{139
  }\cdot\frac{f_2^{811}f_8^2}{f_4^{263}}\nonumber
  \\
  &\qquad +50 \left(
  \frac{1}{f_1^4}\right)^{140}\cdot
   \frac{f_2^{821} }{
  f_4^{265}f_8^2}  \pmod {128}.\label{4-26}
\end{align}
If we
 substitute
\eqref{2-7}
   into \eqref{4-26}
 and exact
 those terms in which the power
  of $q$ is congruent to 1
  modulo 2,
 then divide by $q$ and
 replace
  $q^2
 $ by $q$, we arrive at
\begin{align}\label{4-26-1}
\sum_{n=0}^\infty \overline{bt}(64n+40)q^n
 &\; \equiv  96\frac{f_2^{1683} }{f_1^{1135}
  f_4^{554} }   \pmod {128}.
\end{align}
In light of  \eqref{2-9}, \eqref{4-11}
 and \eqref{4-26-1},
\begin{align}
\overline{bt}(64n+40)\equiv-
 \overline{bt}(16n+10) \pmod {128}.
 \label{4-27}
\end{align}
Replacing $n$ by $9n+2$ in
 \eqref{4-27} and employing
   \eqref{4-16} yields
\begin{align}
\overline{bt}(576n+168)\equiv 0\pmod {128}.
 \label{4-28}
\end{align}
Substituting \eqref{2-6} and \eqref{2-7}
   into \eqref{2-3}
 and exacting
 those terms in which the power
  of $q$ is congruent to 1
  modulo 2,
 then dividing by $q$ and  replacing
  $q^2
 $ by $q$,
  we arrive at
\begin{align}
\sum_{n=0}^\infty \overline{bt}(32n+16)q^n
  \equiv&\;
 64q \left(
  \frac{1}{f_1^4}\right)^{138}
  \cdot
  \frac{1}{f_1^2}\cdot \frac{f_2^{805}
  f_8^2}{ f_4^{259}}
   +112\left(
  \frac{1}{f_1^4}\right)^{139}
  \cdot
  \frac{1}{f_1^2}\cdot\frac{f_2^{815}}{
     f_4^{261}f_8^2}
     \nonumber\\
     &+2\left(
  \frac{1}{f_1^4}\right)^{140}
  \cdot
  \frac{1}{f_1^2}\cdot\frac{f_2^{829}f_8^2}{
       f_4^{275}} \pmod {128}. \label{4-30}
\end{align}
Substituting \eqref{2-6} and \eqref{2-7}
   into \eqref{4-30}
 and picking out
 those terms in which the power
  of $q$ is congruent to 0
  modulo 2,
 then   replacing
  $q^2
 $ by $q$,
  we deduce that
\begin{align}
\sum_{n=0}^\infty \overline{bt}(64n+16)q^n
  \equiv&\;
  64q\frac{f_2^{1675}
  f_8^2}{f_1^{1132}
  f_4^{551}}+64q\frac{
   f_2^{1661}}{f_1^{1128}
    f_4^{537}f_8^2}+114\frac{f_2^{1685}}{
    f_1^{1136}f_4^{553}f_8^2}
    \nonumber\\
    \equiv &\;  114\cdot
    \left(\frac{1}{f_1^4}\right)^{284}\frac{f_2^{1685}}{
     f_4^{553}f_8^2}
    \pmod {128} .
    \qquad ({\rm by}\ \eqref{2-9})\label{4-31}
\end{align}
Substituting   \eqref{2-7}
   into \eqref{4-31}
 and exacting
 those terms in which the power
  of $q$ is congruent to 1
  modulo 2,
 then dividing by $q$ and  replacing
  $q^2
 $ by $q$,
  we arrive at
\begin{align}
\sum_{n=0}^\infty \overline{bt}(128n+80)q^n
    \equiv &\;  96\frac{f_2^{3411}}{f_1^{2287}
     f_4^{1130}}
    \pmod {128} . \label{4-32}
\end{align}
In light of \eqref{2-9},
 \eqref{4-11} and \eqref{4-32},
 \begin{align}
\overline{bt}(128n+80)\equiv -\overline{bt}(16n+10)
 \pmod {128}. \label{4-33}
 \end{align}
 Replacing $n$ by $9n+2$
  in \eqref{4-33}
   and utilizing \eqref{4-16}, we find that for $n\geq 0$,
 \begin{align}
\overline{bt}(1152n+336)\equiv 0
 \pmod {128}.
\label{4-34}
 \end{align}
Replacing $n$ by $72n+21$
 in \eqref{4-1} and using \eqref{4-34}
  yields that for $k\geq 4$ and
  $n\geq 0$,
\begin{align}
\overline{bt}(2^{k}(72n+21) )\equiv 0 \pmod {128}.\label{4-35}
\end{align}
Congruence \eqref{1-4}
 follows from \eqref{v-4}, \eqref{4-16},
  \eqref{4-24},
   \eqref{4-28}  and \eqref{4-35}.

  Substituting
 \eqref{4-15}
  and \eqref{4-14} into \eqref{4-13}
   and extracting those terms in
 which the power of $q$ is
  congruent to 2 mod 3, then dividing by $q^2$ and
 replacing $q^3$ by $q$, we deduce that
\begin{align*}
\sum_{n=0}^\infty \overline{bt}(144n+138)q^n
 &\; \equiv  64 \frac{f_2^6f_6f_{12}}{f_1^3f_4}
 +64 \frac{f_1^7f_6^4}{f_2f_3^2}
\nonumber\\[6pt]
 &\;\equiv 0 \pmod {128},
 \qquad ({\rm by} \eqref{2-9})
\end{align*}
from which, we obtain
\begin{align}
\overline{bt}(144n+138)\equiv 0 \pmod {128}.
 \label{4-36}
\end{align}
Replacing $n$ by $9n+8$ in \eqref{4-23}
 and using \eqref{4-36} yields
\begin{align}\label{4-37}
\overline{bt}(288n+276)\equiv  0\pmod {128}.
\end{align}
Replacing $n$ by $9n+8$ in \eqref{4-27}
 and employing \eqref{4-36} yields
\begin{align}\label{4-38}
\overline{bt}(576n+552)\equiv  0\pmod {128}.
\end{align}
Replacing $n$ by $9n+8$ in \eqref{4-33}
 and utilizing  \eqref{4-36} yields
\begin{align}\label{4-39}
\overline{bt}(1152n+1104)\equiv  0\pmod {128}.
\end{align}
Replacing $n$ by $72n+69$
 in \eqref{4-1} and using \eqref{4-39}
  yields that for $k\geq 4$ and
  $n\geq 0$,
\begin{align}
\overline{bt}(2^{k}(72n+69) )\equiv 0
 \pmod {128}.\label{4-40}
\end{align}
Congruence \eqref{1-5}
 follows from \eqref{v-5} and
 \eqref{4-36}--\eqref{4-38} and \eqref{4-40}.
 This completes the proof. \qed

\section*{Statements and Declarations}

 \noindent{\bf Funding.}
   This work was supported by
 the National Natural Science Foundation of
  China  (grant no.
    12371334) and the Natural Science Foundation of
   Jiangsu Province of China (grant no.
    BK20221383).

 \noindent{\bf Competing Interests.}
 The authors declare that they have
no conflict of interest.

  \noindent{\bf Data Availability Statements.} Data sharing not applicable to this
article as no datasets were generated or analysed during the current
study.


\begin{thebibliography}{99}

\bibitem{Andrews}
G.E. Andrews,
 M.D. Hirschhorn
  and J.A. Sellers,
  Arithmetic properties of partitions
   with even parts
distinct, Ramanujan J.
 23 (2010)
  169--181.

\bibitem{Berndt-1}
B.C. Berndt, Ramanujan's Notebooks, Part III, Springer, New York,
1991.

\bibitem{Chan-1}
    H.-C. Chan, Ramanujan's cubic continued
     fraction and an analog of his
      ``most beautiful identity",
       Int. J. Number Theory 6 (2010)
        673--680.

   \bibitem{Chan-2}
     H.-C. Chan, Ramanujan's
   cubic continued fraction and Ramanujan type congruences for a
certain partition function, Int. J. Number Theory
 6 (2010) 819--834.

\bibitem{Chan-3} H.-C. Chan,
Distribution of a certain partition function modulo powers of
primes, Acta Math. Sin. (Engl. Ser.) 27 (2011)
 625--634.

\bibitem{Hirschhorn}
M.D. Hirschhorn  and J.A. Sellers,
 A congruence modulo 3 for
partitions into distinct non-multiples
 of four, J. Integer
Sequences 17 (2014)  Article 14.9.6.




\bibitem{Kim-2010}
B. Kim, The overcubic partition function mod 3, Ramanujan
Rediscovered, Ramanujan Math.
 Soc. Lect. Notes Ser. 14 (2010)
157--163.

\bibitem{Kim-2012}
B. Kim, On partition congruences
 for overcubic partition pairs,
  Commun. Korean Math. Soc.
   27 (2012)  477--482.

\bibitem{Saikia}
 M. P. Saikia  and A. Sarma,
Further arithmetic properties of overcubic
partition triples,
    Bull. Aust. Math. Soc.
     to appear (https://doi.org/10.1017/S000497272400114X).

\bibitem{Sellers}
J.A. Sellers, Elementary proofs of congruences for the cubic and
overcubic partition functions,
Austral. J. Combinatorics 60 (2)
(2014)  191--197.

\bibitem{Nayaka}
S.S. Nayaka,
 B.N. Dharmendra
 and M.C.M. Kumar, Divisibility properties for
overcubic partition triples, Integers
 24 (2024) \# A80.

\bibitem{Zhao}
H. Zhao and Z.
 Zhong,    Ramanujan type congruences for a
partition function,
  Electron. J. Combin. 18  (2011) \# P58.


\end{thebibliography}
\end{document}